\documentclass[11pt,a4paper,reqno]{amsart}
\usepackage{graphicx}
\usepackage{epsfig}
 \theoremstyle{definition}
 
 \theoremstyle{remark}
 
\begin{document}
\setcounter{page}{1}
\begin{flushleft}
{\scriptsize 
}
\end{flushleft}
\bigskip
\bigskip
\title[S. S. Siddiqi, M. Iftikhar: Variational Iteration Homotopy ... ] {Variational Iteration Homotopy Perturbation Method
 for the Solution
of Seventh Order Boundary Value Problems }
\author[]{Shahid S. Siddiqi$^1$, Muzammal Iftikhar$^2$}
\thanks{$^{1,2}$Department of Mathematics, University of the Punjab, Lahore 54590,
 Pakistan,\\ \indent\,\,\, $^2$Department of Mathematics, University of Education, Okara
Campus, Okara 56300,
 Pakistan
\\ \indent\,\,\,$^1$e-mail: shahidsiddiqiprof@yahoo.co.uk,\\ \indent\,\,\,$^2$e-mail:
miftikhar@hotmail.com.
}
\begin{abstract}
The induction motor behaviour is represented by a fifth order
differential equation model. Addition of a torque  correction factor
to the model accurately reproduces the transient torques and
instantaneous real and reactive power flows of the full seventh
order differential equation model. The variational iteration
homotopy perturbation method (VIHPM) is employed to solve the
seventh order boundary value problems. The approximate solutions of
the problems are obtained in terms of a rapidly convergent series.
Several numerical examples are given to illustrate the
implementation and efficiency of the method.

\bigskip
\noindent Keywords: Variational iteration homotopy perturbation
method; Boundary value problems; Linear and nonlinear problems;
Approximate solution.

\bigskip \noindent AMS Subject Classification: 34B05; 34B10; 34B15.

\end{abstract}
\maketitle

\smallskip
\section*{Introduction}The theory of seventh order boundary value problems
is not much available in the numerical analysis literature.  These
problems  are generally arise in modelling induction motors with two
rotor circuits.\\  \indent The induction motor behavior is
represented by a fifth  order differential equation model. This
model contains two stator state variables, two rotor state variables
and one shaft speed. Normally, two more variables must be added to
account for the effects of a second rotor circuit representing deep
bars, a starting cage or rotor distributed parameters. \\  \indent
To avoid the computational burden of additional state variables when
additional rotor circuits are required, model is often limited to
the fifth order and rotor impedance is algebraically altered as
function of rotor speed under the assumption that the frequency of
rotor currents depends on rotor speed. This approach is efficient
for the steady state response with sinusoidal voltage, but it does
not hold up during the transient conditions, when rotor frequency is
not a single value. So, the behavior of such models show up in the
seventh order (Richards and Sarma, 1994)
\\  \indent Siddiqi and Ghazala (Siddiqi and Ghazala, 2006a and 2006b)
 ~presented the solutions of fifth and sixth order boundary value
problems using non-polynomial spline. Noor and Mohyud-Din (Noor and
Mohyud-Din, 2009)
  ~applied modified variational iteration
method for solving the boundary layer problem in unbounded domain.
Matinfar $\textit{et~al}$ (Matinfar $\textit{et~al}$, 2010) 
implemented the variational homotopy perturbation method to obtain
the solution of Fisher's equation. (Siddiqi and Iftikhar, 2012) used
the variation of parameter method to solve the seventh order
boundary value problems. \\ \indent Odibat  discussed the
convergence of variational iteration method in (Odibat, 2010).
Tatari and Dehghan presented the sufficient conditions to guarantee
the convergence of the variation iteration method (Tatari and
Dehghan, 2007).
\\  \indent
 The aim of this study is to solve the seventh order
boundary value problems and the variational iteration homotopy
perturbation method is used for this purpose.

\bigskip
\section*{ Variational Iteration Homotopy Perturbation
Method (Noor and
Mohyud-Din, 2009,  Matinfar $\textit{et~al}$, 2010, Mohyud-Din $\textit{et~al}$, 2010)  }

Variational iteration homotopy perturbation method is formulated by
the coupling of variational iteration method and homotopy
perturbation method. The boundary value problem is considered as
under
\begin{eqnarray}\label{e2.1}
 L[u(x)]+N[u(x)]= g(x),
\end{eqnarray}
where $L$ and $N$ are linear and nonlinear operators respectively
and $g(x)$ is a forcing term. Following the variational iteration
method used by He (He, 1998, 1999, 2000a, 2001). 
The correct functional for the problem (\ref{e2.1}) can be written
as
\begin{equation}\label{e2.2}
 u_{n+1}(x)=u_{n}(x)+\int_0^x
 \lambda(t)\{Lu_{n}(t))+N\widetilde{u}_{n}(t)-g(t)\}dt,
\end{equation}
where $\lambda$ is a Lagrange multiplier, that can be identified
optimally $via$ variational iteration method.
 Here, $\widetilde{u}_{n}$ is considered to be a restricted variation
which shows that $\delta\widetilde{u}_{n}=0$. Making the correct
functional (\ref{e2.2}) stationary, yields
\begin{eqnarray}\label{ea2.6.2}
 \delta v_{n+1}(x)&=&\delta v_{n}(x)+\delta \int_0^x
 \lambda(t)\{Lv_{n}(t))+N\widetilde{v}_{n}(t)-g(t)\}dt\nonumber\\
 &=& \delta v_{n}(x)+ \int_0^x
 \delta \{\lambda (t)Lv_{n}(t))\}dt.
\end{eqnarray}
Its stationary conditions can be obtained using integration by parts
in Eq. (\ref{ea2.6.2}). Therefore, the Lagrange multiplier can be
 written as
  \begin{equation}\label{ea2.6.6}
 \lambda=
 \frac{(-1)^{m}(t-x)^{m-1} }{(m-1)!}.
\end{equation} Applying the homotopy
perturbation method the following relation is obtained as
\begin{equation}\label{e2.3}
\sum^{\infty}_{i=0}p^{i}u_{i}(x)=u_{0}(x)+\int_0^x
 \lambda(t)\left\{L \left(\sum^{\infty}_{i=0}p^{i}u_{i}\right)+N(\sum^{\infty}_{i=0}p^{i}\widetilde{u}_{i})\right\}dt-\int_0^x\lambda(t)g(t)dt,
\end{equation}
 Equating the like powers of  $p$ gives
$u_{0},~u_{1},~\cdots$. The embedding parameter $p\in[0, 1]$ can be
used as an expanding parameter. The approximate solution of the
problem \eqref{e2.1}, therefore, can be expressed as
\begin{equation}\label{e2.4}
 u=\lim_{p\rightarrow1} \sum^{\infty}_{i=0}p^{i}u_{i} =u_{0}+u_{1} + u_{2}+\cdots
\end{equation}
The series \eqref{e2.4} is convergent for most of the cases
It is assumed that \eqref{e2.4} has a unique solution.\\In fact, the
solution of the problem (\ref{e2.1}) is considered as the fixed
point of the following functional under the suitable choice of the
initial term $v_0(x)$.\begin{equation}\label{e2.6.2b}
 v_{n+1}(x)=v_{n}(x)+\int_0^x
 \lambda(t)\{\mathbf{L}v_{n}(t))+\mathbf{N}{v}_{n}(t)-g(t)\}dt.
\end{equation}
\section*{ Convergence }
In this section, we will present Banach's theorem about the
convergence of the VIHPM. The VIHPM changes the given differential
equation into a recurrence sequence of functions. The limit of this
sequence is
considered as the solution of the given differential equation\\
\textbf{Theorem 1}(Banach's fixed point theorem) (Tatari and
Dehghan, 2007)~~Suppose that $X$ is a Banach space and $B : X
\longrightarrow X$ is a nonlinear mapping, and assume that
\begin{eqnarray}\label{ea}
\|B[u]-B[\bar{u}]\|\leq \gamma \|u-\bar{u} \|, ~~ \forall
~u,\bar{u}\in X.
\end{eqnarray}
for some constant $\gamma <1$. Then $B$ has a unique fixed point.
Moreover, the sequence
\begin{eqnarray}\label{eb}
u_{n+1}=B[u_{n}]
\end{eqnarray}
with an arbitrary choice of $u_0 \in X$ converges to the fixed point
of $B$ and
\begin{eqnarray}\label{ec}
 \|u_k-u_1\|\leq  \sum _{j=l-1}^{k-2}\gamma ^{j} \|u_1-u_0
\|,
\end{eqnarray}
According to Theorem 1, for the nonlinear mapping
\begin{equation}\label{ec}
 B[u]=u(x)+\int_0^x
 \lambda(t)\{Lu_{n}(t))+N\widetilde{u}_{n}(t)-g(t)\}dt,
\end{equation}
is a sufficient condition for convergence of the variational
iteration homotopy perturbation method is strictly contraction of
$B$. Furthermore, the sequence (\ref{eb}) converges to the fixed
point of $B$ which also is the solution of the problem (1).\\
To implement the method, some numerical examples are considered in
the following section.
\bigskip
\section*{ Numerical Examples }
\textbf{Example 1}~~The following seventh order linear boundary
value problem is considered
\begin{eqnarray}\label{e4.1}
\left.\begin{array}{lll}
 u^{(7)}(x)=-u(x)-e^{x}(35+12x+2x^2), 0 \leq  x \leq1,\\
u(0)=0,~~~~~~~~~~~~~~~u(1)=0,~~~\\
u^{(1)}(0)=1,~~~~~~~~~~~~u^{(1)}(1)=-e,~~~~~~~~~~~~~~~~~~~~~~\\
u^{(2)}(0)=0,~~~~~~~~~~~~u^{(2)}(1)=-4e,\\
u^{(3)}(0)=-3.~~~~~~~~~~~~~~~~~~~~~~~~~~~~~~\\
\end{array} \right\}\
\end{eqnarray}
The exact solution of the Example 1  is
 $u(x)=x(1-x)e^x,$ (Siddiqi and Iftikhar, 2013).\\
 The correct functional for the
problem (\ref{e4.1}) can be written as
\begin{equation}\label{e2.6.2aa}
 v_{n+1}(x)=v_{n}(x)+\int_0^x
 \lambda(t)\{v^{(7)}_{n}(t))-v_{n}(t)+e^{t}(35+12t+2t^2)\}dt,
\end{equation}
Making the correct functional (\ref{e2.6.2aa}) stationary, yields
\begin{eqnarray}\label{e2.6.2bb}
 \delta v_{n+1}(x)&=&\delta v_{n}(x)+\delta \int_0^x
 \lambda(t)\{v^{(7)}_{n}(t))-v_{n}(t)+e^{t}(35+12t+2t^2)\}dt\nonumber\\
 &=&\delta v_{n}(x)+ \int_0^x \delta
 \{ \lambda(t)v^{(7)}_{n}(t))t\}dt,
\end{eqnarray}
Hence, the following stationary conditions can be
determined\begin{eqnarray*}\label{e12.6.6}
  \lambda^{(7)}(t) &=&0,\\\lambda (t)|_{t=x} &=& 0,\\ \lambda' (t)|_{t=x} &=&  0,\\ 
&\vdots & \\  \lambda^{(5)} (t)|_{t=x} &=&  0,\\
~1+\lambda^{(6)}(t)|_{t=x} &=&  0,
\end{eqnarray*}
which yields
   \begin{equation}\label{e2.6.6}
 \lambda=
 \frac{(-1)^{7}(t-x)^{6} }{(6)!}.
\end{equation}
The Lagrange multiplier can be identified as follows
 \begin{equation}\label{e4.2}
 \lambda=
 \frac{(-1)^{m}(t-x)^{m-1} }{(m-1)!}.
\end{equation}
 According to (\ref{e2.3}), the following iteration formulation is
obtained
\begin{equation}\label{e4.3}
\sum^{\infty}_{i=0}p^{i}u_{i}(x)=u_{0}(x)+\int_0^x
\frac{(-1)^{7}(t-x)^{6} }{6!}
\left\{\left(\sum^{\infty}_{i=0}p^{i}u_{i}\right)^{(7)}+\sum^{\infty}_{i=0}p^{i}u_{i}+e^{t}(35+12t+2t^2))\right\}dt.
\end{equation}
Now, assume that an initial approximation has the form
\begin{equation}\label{e4.4}
 u_{0}(x)=x-\frac{x^3}{2}+Ax^4+Bx^5+Cx^6.
\end{equation}
Comparing the coefficient of like powers of $p$
\begin{eqnarray*}
p^{0}:\quad u_{0}(x)&=&x-\frac{x^3}{2}+Ax^4+Bx^5+Cx^6, \nonumber \\
p^{1}:\quad u_{1}(x)&=&-\frac{x^7}{144}-\frac{x^8}{840}-\frac{x^9}{5760}-\frac{x^{10}}{45360}
+\left(\frac{-107}{39916800}-\frac{A}{1663200}\right)x^{11}\nonumber\\
&&+\left(\frac{-1}{3548160}-\frac{B}{3991680}\right)x^{12}+O(x)^{13},\\
\vdots
\end{eqnarray*}
 where $A$, $B$ and
$C$ are unknown constants to be
 determined later.\\
Using the first two approximations the series solution cane be
written as
\begin{eqnarray*}
u(x)&=&
x-\frac{x^3}{2}+Ax^4+Bx^5+Cx^6-\frac{x^7}{144}-\frac{x^8}{840}-\frac{x^9}{5760}-\frac{x^{10}}{45360}
+\left(\frac{-107}{39916800}-\frac{A}{1663200}\right)x^{11}\nonumber\\
&&+\left(\frac{-1}{3548160}-\frac{B}{3991680}\right)x^{12}+O(x)^{13}.
\end{eqnarray*}
Using the boundary conditions (\ref{e4.1}),  the values of the
unknown constants can be determined as follows \\
$A=-0.3333333170467781$,~~~~~$B=-0.12500003614813987$,~~~~~
$C=-0.03333331303032349.$\\Finally, the series solution is
\begin{eqnarray*}
u(x)&=&x -(0.5)x^3 -0.333333x^4 -0.125x^5 -0.0333333x^6
-0.00694444x^7 \\&&- 0.00119048x^8- 0.000173611x^9
-0.0000220459x^{10} - \left(2.48016\times10^{-6}\right)x^{11}
\\&&- \left(2.50521\times10^{-7}\right)x^{12}
+O(x)^{13}.
\end{eqnarray*}
The comparison of the values of maximum absolute errors of the
present method with the variation of parameter method (Siddiqi and
Iftikhar, 2013) for the Example 1 is given in Table 1, which shows
that the present method is more accurate. In Figure 1 the comparison
of exact and approximate solutions
is given and absolute errors are plotted in Figure 2 for Example 1.\\\\
\textbf{Example 2}~~The following seventh order nonlinear boundary
value problem is considered
\begin{eqnarray}\label{e4.6}
\left.\begin{array}{lll} u^{(7)}(x) =e^{-x}u^{2}(x),~~ 0 <
x<1,~~~~~~~~~~~~~~~~~~~~~~~~~~~~~~~~~~~~~~~~~~~~~~~\\
u(0)=u^{(1)}(0)=u^{(2)}(0)=u^{(3)}(0)=1,~~~~~~~~~~~~~~~~~~~~~~~~~~~~~~~~~~~~\\
u(1)=u^{(1)}(1)=u^{(2)}(1)=e.~~~~~~~~~~~~~~~~~~~~~~~~~~~~~~~~~~~~~~~~~~~~~~~~~\\
\end{array} \right\}\
\end{eqnarray}
The exact solution of the Example 2 is
 $u(x)=e^x,$ (Siddiqi and Iftikhar, 2013)\\
 The Lagrange multiplier can be identified as follows
 \begin{equation}\label{e4.7}
 \lambda=
 \frac{(-1)^{m}(t-x)^{m-1} }{(m-1)!}.
\end{equation}
 According to (\ref{e2.3}), the following iteration formulation is
obtained
\begin{equation}\label{e4.8}
\sum^{\infty}_{i=0}p^{i}u_{i}(x)=u_{0}(x)+\int_0^x
\frac{(-1)^{7}(t-x)^{6}}{6!}
\left\{\left(\sum^{\infty}_{i=0}p^{i}u_{i}\right)^{(7)}-\left(\sum^{\infty}_{i=0}p^{i}u_{i}\right)^2e^{-t}\right\}dt.
\end{equation}
Now, assume that an initial approximation has the form
\begin{equation}\label{e4.9}
 u_{0}(x)=1+x+\frac{x^2}{2}+\frac{x^3}{6}+Ax^4+Bx^5+Cx^6.
\end{equation}
Comparing the coefficient of like powers of $p$
\begin{eqnarray*}
p^{0}:\quad u_{0}(x)&=&1+x+\frac{x^2}{2}+\frac{x^3}{6}+Ax^4+Bx^5+Cx^6, \nonumber \\
p^{1}:\quad u_{1}(x)&=&\frac{x^7}{5040}+\frac{x^8}{40320}+\frac{x^9}{362880}+\frac{x^{10}}{3628800}+\left(\frac{-1}{39916800}+\frac{A}{831600}\right)x^{11}\nonumber\\
&&+\left(\frac{-1}{479001600}+\frac{B}{1995840}\right)x^{12}+O(x)^{13},\\
\vdots
\end{eqnarray*}
 where $A$, $B$ and
$C$ are unknown constants to be
 determined later.\\
Using the first two approximations the series solution cane be
written as
\begin{eqnarray*}
u(x)&=&
1+x+\frac{x^2}{2}+\frac{x^3}{6}+Ax^4+Bx^5+Cx^6+\frac{x^7}{5040}+\frac{x^8}{40320}+\frac{x^9}{362880}+\frac{x^{10}}{3628800}\nonumber\\
&&+\left(\frac{-1}{39916800}+\frac{A}{831600}\right)x^{11}+\left(\frac{-1}{479001600}+\frac{B}{1995840}\right)x^{12}+O(x)^{13}.
\end{eqnarray*}
Using the boundary conditions (\ref{e4.6}),  the values of the
unknown constants can be determined as follows \\
$A=0.041666667529862395$,~~~~~$B=0.008333331197193119$,~~~~~
$C=0.001388890268167299.$\\ Finally, the series solution is
\begin{eqnarray*}
u(x)&=&1. +  x + (0.5)x^2 + 0.166667x^3 + 0.0416667x^4 +
0.00833333x^5 +0.00138889x^6 \\&&+ 0.000198413x^7 + 0.0000248016x^8
+(2.75573\times10^{-6})x^9 + (2.75573\times10^{-7})x^{10}
\\&&+ (2.50521\times10^{-8})x^{11} + (2.08768\times10^{-9})x^{12}+
O(x)^{13}.
\end{eqnarray*} In Table 2,  errors obtained by the
present method are compared with errors obtained using the variation
of parameters method (Siddiqi and Iftikhar, 2013) for the Example 2.
It is observed that the maximum absolute error value for the present
method is $4.5614\times10^{-9}$ which is better than the maximum
absolute error value , $7.7176\times10^{-7}$, of the variation of
parameters method (Siddiqi and Iftikhar, 2013). Figure 2 shows the
comparison of exact and approximate solutions and absolute errors
are plotted in Figure 4 for Example 2.
 The results reveals that the present method is more accurate.\\
\textbf{Example 3}~~The following seventh order nonlinear boundary
value problem is considered
\begin{eqnarray}\label{e4.7}
\left.\begin{array}{lll}
 u^{(7)}(x)&=&-u(x)u'(x)+g(x),~~0 \leq  x \leq1,~~~~~~~~~~~~~~~~~~~~~~~~\\
u(0)&=&0,~~~~u(1)=0,~~~~~~~~~~~~~~~~~~~~~~~~~~~~~~~~~~~~~~~~~~~~~\\
u^{(1)}(0)&=&1,~~~u^{(1)}(1)=-e,~~~~~~~~~~~~~~~~~~~~~~~~~~~~~~~~~~~~~~~~~\\
u^{(2)}(0)&=&0,~~~~~~u^{(2)}(1)=-4e,~~~~~~~~~~~~~~~~~~~~~~~~~~~~~~~~~~~~\\
u^{(3)}(0)&=&-3.~~~~~~~~~~~~~~~~~~~~~~~~~~~~~~~~~~~~~~~~~~~~~~~~~~~~~~~~~~~~~~\\
\end{array} \right\}\
\end{eqnarray}
where $g(x)=e^x(-35+(-13+e^x)x-(1+2e^x)x^2+e^x x^4)$.\\ The exact
solution of the Example 3.3  is
  $u(x)=x(1-x)e^x.$\\
  Following the procedure of the previous examples the series solution, using the first two
  approximations, can be written as
\begin{eqnarray*}
u(x)&=&x -(0.5)x^3 -0.333333x^4 -0.125x^5 -0.0333332x^6
-0.00694444x^7 \\&&- 0.00119048x^8 -0.000173611x^9
-0.0000220459x^{10} - (2.48016\times10^{-6})x^{11}
\\&&- (2.50521\times10^{-7})x^{12}
+O(x)^{13}.
\end{eqnarray*}
The comparison of the exact solution with the series solution of the
Example 3 is given in Table 3. In Figure 5 absolute errors are
plotted.\\\\
\textbf{Example 4}~~The following seventh order nonlinear three
point boundary value problem is considered
\begin{eqnarray}\label{e4.8}
\left.\begin{array}{lll}
 u^{(7)}(x)&=&u(x)u'(x)-e^x (6+x-e^xx+e^xx^2),~~~~ 0 \leq  x \leq1,\\
u(0)&=&1,~~~~u(\frac{1}{2})=\frac{e^{\frac{1}{2}}}{2},~~~~~~~~~~~~~~~~~~~~~~~~~~~~~~~~~~~~~~~~~~~~~\\
u^{(1)}(0)&=&0,~~~u^{(1)}(\frac{1}{2})=-\frac{e^{\frac{1}{2}}}{2},~~~~~~~~~~~~~~~~~~~~~~~~~~~~~~~~~~~~~~~~\\
u^{(2)}(0)&=&-1,~~~u^{(2)}(1)=-2e,~~~~~~~~~~~~~~~~~~~~~~~~~~~~~~~~~~~~~~~~~\\
u(1)&=&0.~~~~~~~~~~~~~~~~~~~~~~~~~~~~~~~~~~~~~~~~~~~~~~~~~~~~~~~~~~~~~~~~\\
\end{array} \right\}\
\end{eqnarray}
The exact solution of the Example 4  is
 $u(x)=(1-x)e^x.$\\
 Following previous examples the series solution, using the first two
  approximations, can be written as
\begin{eqnarray*}
u(x)&=&1. -(0.5)x^2 -0.333333x^3 -0.125001x^4 -0.0333319x^5
-0.00694445x^6 \\&&- 0.00119048x^7 -0.000173611x^8 -0.0000220459x^9
- (2.48016\times10^{-6})x^{10}
\\&&- (2.50517\times10^{-7})x^{11} - (2.29651\times10^{-8})x^{12}
+O(x)^{13}.
\end{eqnarray*}
The comparison of the exact solution with the series solution of the
Example 4 is given in Table 4. Absolute errors are plotted
in Figure 6.\\
\hspace{1mm}
\begin{center}
\includegraphics[width=6cm]{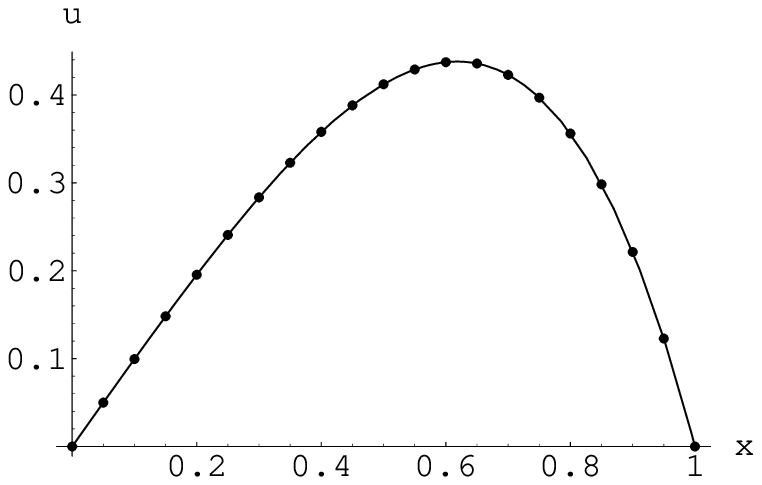}\\
{\footnotesize Figure 1. Comparison between the exact solution and
the approximate solution for Example  (1). Dotted line: approximate
solution, solid line: the exact solution.}
\end{center}\label{figa} \hspace{3mm}

\hspace{1mm}
\begin{center}
\includegraphics[width=6cm]{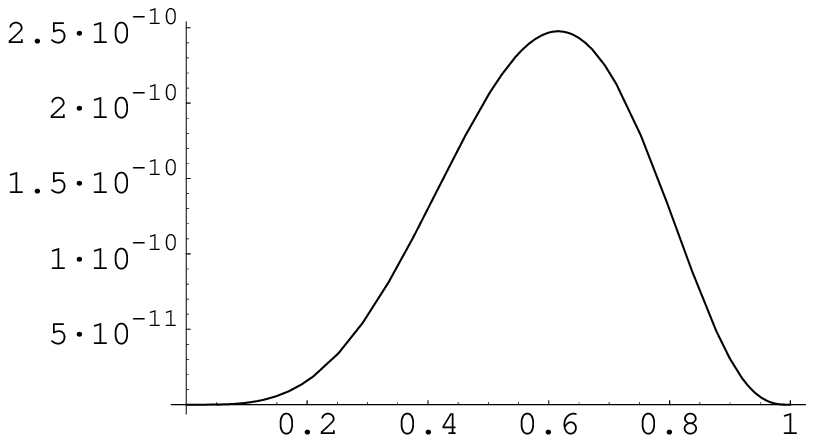}\\
{\footnotesize Figure 2. Absolute errors for Example (1).}
\end{center}\label{figb)} \hspace{3mm}

\hspace{1mm}
\begin{center}
\includegraphics[width=6cm]{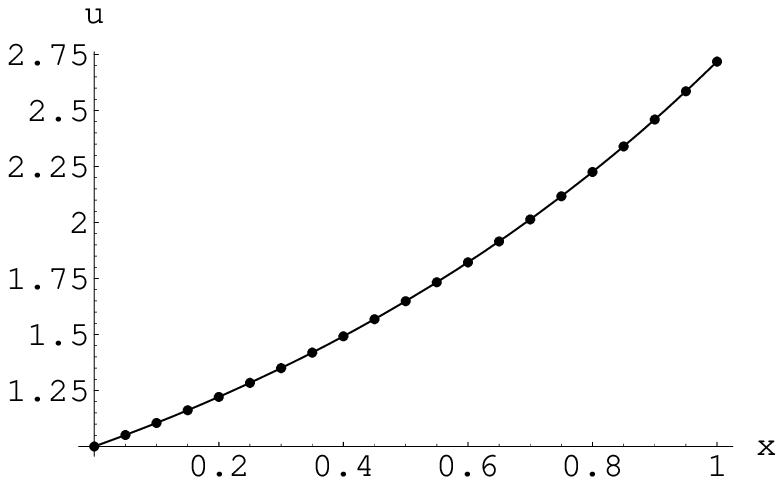}\\
{\footnotesize Figure 3. Comparison between the exact solution and
the approximate solution for Example (2). Dotted line: approximate
solution, solid line: the exact solution.}
\end{center}\label{figc} \hspace{3mm}

\hspace{1mm}
\begin{center}
\includegraphics[width=6cm]{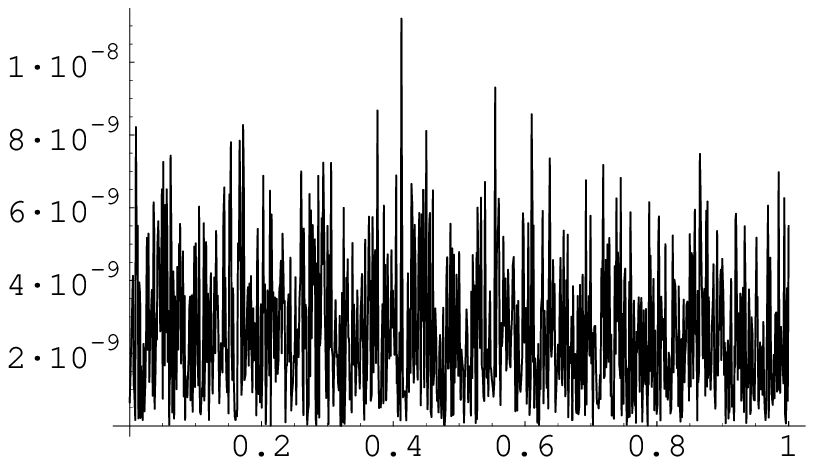}\\
{\footnotesize Figure 4. Absolute errors for Example  (2).}
\end{center}\label{figd} \hspace{3mm}

\hspace{1mm}
\begin{center}
\includegraphics[width=6cm]{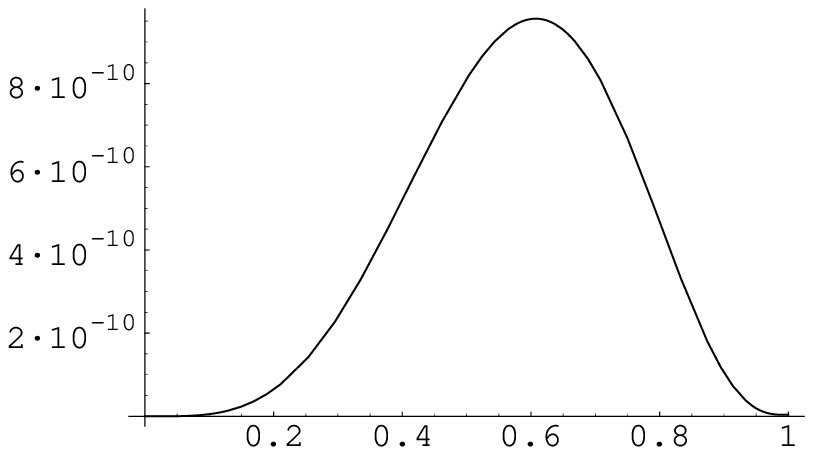}\\
{\footnotesize Figure 5. Absolute errors for Example  (3).}
\end{center}\label{fige} \hspace{3mm}

\hspace{1mm}
\begin{center}
\includegraphics[width=6cm]{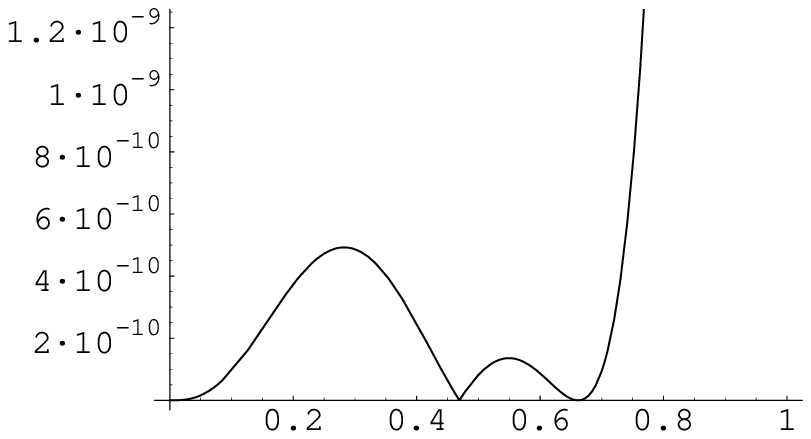}\\
{\footnotesize Figure 6. Absolute errors for Example  (4).}
\end{center}\label{figf} \hspace{3mm}

\bigskip
\section{Conclusion}

In this paper, variational iteration homotopy perturbation method
has been applied to obtain the numerical solutions of linear and
nonlinear seventh order boundary value problems. The variational
iteration homotopy perturbation method solves nonlinear problems
without using He's or Adomian's polynomials. The method gives
rapidly converging series solutions in both linear and nonlinear
cases. The numerical results revealed that the present method is a
powerful mathematical tool for the solution of seventh order
boundary value problems. Numerical examples also show the accuracy
of the method.

\begin{table}[bht]
\caption{Comparison of maximum absolute errors for Example 1}
\centering
\begin{small}
\begin{tabular}{|c|c|}
\hline  present method& Variation of Paramters method \\
&(Siddiqi and Iftikhar, 2013)\\
\hline$2.46782\times 10^{-10}$&$2.1729\times 10^{-09}$\\
\hline
\end{tabular}\end{small}
\end{table}
\begin{table}[bht]
\caption{Comparison of numerical results for Example 2} \centering
\begin{small}
\begin{tabular}{|c|c|c|c|c|}
\hline $x$& Exact solution  & Approximate & Absolute Error&Absolute
Error
\\&& solution&present method&(Siddiqi and Iftikhar, 2013)\\
\hline$0.0$&1.0000&1.0000&1.32455E-09&0.0000\\
\hline$0.1$&1.1051&1.1051&5.26137E-10&2.26257E-07\\
\hline$0.2$&1.2214& 1.2214&1.64015E-09&4.38942E-07\\
\hline$0.3$&1.3498& 1.3498&4.56139E-09&6.1274E-07\\
\hline$0.4$&1.4918& 1.4918&2.9619E-09&7.71759E-07\\
\hline$0.5$&1.6487& 1.6487&7.54889E-10&7.71759E-07\\
\hline$0.6$&1.8221& 1.8221&2.67612E-09&7.37682E-07\\
\hline$0.7$&2.0137& 2.0137&8.42306E-10&6.25932E-07\\
\hline$0.8$&2.2255& 2.2255&1.16866E-09&4.68244E-07\\
\hline$0.9$&2.4596&2.4596&4.4716E-09&2.95852E-07\\
\hline$1.0$&2.7183&2.7183&1.02746E-09&1.25922E-07\\
\hline
\end{tabular}\end{small}
\end{table}
\begin{table}[bht] \caption{Comparison of numerical results for
Example 3} \centering
\begin{small}
\begin{tabular}{|c|c|c|c|}
\hline $x$& Exact solution  & Approximate & Absolute Error
\\&&solution&\\
\hline$0.0$&0.0000&0.0000&0.00000\\
\hline$0.1$&0.0994654&0.0994654&5.28944E-12\\
\hline$0.2$&0.195424& 0.195424&6.44606E-11\\
\hline$0.3$&0.28347& 0.28347&2.38427E-10\\
\hline$0.4$&0.358038& 0.358038&5.20559E-10\\
\hline$0.5$&0.41218& 0.41218&8.11431E-10\\
\hline$0.6$&0.437309&0.437309&9.55209E-10\\
\hline$0.7$&0.422888& 0.422888&8.30543E-10\\
\hline$0.8$&0.356087& 0.356087&4.67351E-10\\
\hline$0.9$&0.221364&0.221364&1.04882E-10\\
\hline$1.0$&0.0000&3.90259E-12&3.90259E-12\\
\hline
\end{tabular}\end{small}
\end{table}
\begin{table}[bht]
\caption{Comparison of numerical results for Example 4} \centering
\begin{small}
\begin{tabular}{|c|c|c|c|}
\hline $x$& Exact solution  & Approximate  & Absolute Error \\
&&series solution&\\
\hline$0.0$&1.0000&1.0000&0.0000\\
\hline$0.1$&0.9946&0.9946&9.48615E-11\\
\hline$0.2$&0.9771& 0.9771&3.7371E-10\\
\hline$0.3$&0.9449& 0.9449&4.8626E-10\\
\hline$0.4$&0.8950& 0.8950&2.46565E-10\\
\hline$0.5$&0.8243& 0.8243&8.16711E-11\\
\hline$0.6$&0.7288& 0.7288&8.67514E-11\\
\hline$0.7$&0.6041& 0.6041&9.51461E-11\\
\hline$0.8$&0.4451& 0.4451&2.63398E-09\\
\hline$0.9$&0.2459&0.2459&1.44494E-08\\
\hline$1.0$&0.0000&-4.90417E-08&4.90417E-08\\
\hline
\end{tabular}\end{small}
\end{table}

\end{document}